\documentclass{l4dc2026}

\usepackage[utf8]{inputenc} 
\usepackage[T1]{fontenc}    
\usepackage{hyperref}       
\usepackage{url}            
\usepackage{booktabs}       
\usepackage{nicefrac}       
\usepackage{microtype}      
\usepackage{lipsum}
\usepackage{graphicx}
\usepackage{tikz}
\usepackage{layouts}
\usepackage{pgfplots}
\pgfplotsset{compat=1.14}

\usepackage{mathtools}
\usetikzlibrary {arrows.meta}

\title{Learning Nonholonomic Dynamics with Constraint Discovery}
\usepackage{times}

\author{
 \Name{Baiyue Wang}
  \Email{baiyuew@umich.edu} \\
  \addr Ann Arbor, MI 48104
   \AND
 \Name{Anthony Bloch} 
  \Email{abloch@umich.edu} \\
  \addr Department of Mathmatics, University of Michigan, Ann Arbor, MI 48109 
}

\begin{document}
\maketitle           
\begin{abstract}
We consider learning nonholonomic dynamical systems while discovering the constraints, and describe in detail the case of the rolling disk. A nonholonomic system is a system subject to nonholonomic constraints. Unlike holonomic constraints,  nonholonomic constraints do not define a sub-manifold on the configuration space. Therefore, the inverse problem of finding the constraints has to involve the tangent bundle. This paper discusses a general procedure to learn the dynamics of a nonholonomic system through Hamel's formalism, while discovering the system constraint by parameterizing it, given the data set of discrete trajectories on the tangent bundle $TQ$. We prove that there is a local minimum for convergence of the network. We also preserve symmetry of the system by reducing the Lagrangian to the Lie algebra of the selected group.
\end{abstract}


\section{Introduction}

The theory of symmetries has been generalized to describe the invariant quantities in a dynamical system under transformations. It has caught the attention of researchers in the field of machine learning due to the increasing need for understanding the structure of neural networks. Although it seems that neural networks are far from symmetric, there can be some innate symmetries \cite{bloch2024symmetric}, and a network can be designed to be symmetric \cite{ijrr2024_huang, varghese2024sympgnns}.

Some research has studied how to discover the symmetries of a system through machine learning. Most of the them are based on symplecticity and Hamiltonian systems. One example of learning the dynamics is  learning the evolution of the system: learning the symplectic map \cite{pmlr-chen21r}, learning the symplectic form \cite{neurips_chen21}, learning through the Lie-Poisson brackets \cite{eldred2023liepoisson, eldred2024clpnets, gruber2023reversible}, learning with a symplectic network design \cite{varghese2024sympgnns}, or learning through the $G$-invariant Lagrangian submanifold \cite{vaquero2023symmetry}. These works were conducted based on the assumption that the system is symmetric. While we can make that assumption, there are ways to discover the symmetry automatically through an equivariant network design such as a Lie algebra convolutional network \cite{neurips2021_dehmamy}, or dictionary learning \cite{ghosh2023dictionarylearning}.

A rather less studied topic is to learn a  symmetric system with constraints. Constraints can be holonomic (defining a submanifold 
of the configuration space) or nonholonomic ( constraint on the velocities which are not integrable). The difficulties of this problem lie in that one needs to learn the constraints implicitly from the dynamics of the system. Learning a holonomic constraint (a constraint on the system configurations) can be seen as learning the null space of the constraint equation \cite{icra2015_lin}, in which case the system symmetry is not considered. 

In this work, we propose a general method for learning nonholonomic systems, and at the same time not ignoring the system symmetry, by leveraging Hamel's formalism of nonholonomic dynamics with symmetry \cite{Bloch2009Quasi}. Hamel's formalism treats the constraints intrinsically  as a distribution in the tangent bundle rather than imposing the constraints via  multipliers; therefore, the dynamics is interpreted as an ordinary differential equation without an additional constraint. In this formalism we train the network in a similar manner to a Neural ODE \cite{chen2019neuralode}.

\section{Problem Statement}
Let us first define a nonholonomic constraint.
\begin{definition}[\cite{Bloch2015}]
    Consider a mechanical system, subject to  linear velocity constraints that in generalized coordinates can be expressed as 
    \begin{equation}
    \label{eq:def-constraint}
        a(q)\dot{q} = [a_1(q), \dots, a_n(q)]\dot{q}=0,
    \end{equation}
    where $\dot{q}$ is regarded as a column vector, and $a(q)$ is a matrix of $n$ columns. A constraint is said to be \textbf{nonholonomic} or \textbf{non-integrable} if there is no function $h$ of $q$ such that the constraint can be written as $h(q) = \mathrm{constant}$, or $(\partial h / \partial q^i)\dot{q}^i = 0$.
\label{def:nhc-constraint}
\end{definition}

In the rolling disk dynamics, which will be introduced in more detail in Section \ref{sec:rolling-penny}, the constraints have  the form 
\begin{equation}
\label{eq:constraint}
    \begin{bmatrix}
        -R\cos{\varphi} & 0 & 1 & 0 \\
        -R\sin{\varphi} & 0 & 0 & 1
    \end{bmatrix} \begin{bmatrix}
        \dot{\theta} \\ \dot{\varphi} \\ \dot{x} \\ \dot{y}
    \end{bmatrix} = 0.
\end{equation}

One can see that the constraints are non-integrable.

Suppose that we are given a set of trajectory data of the nonholonomic system $(q^i, \dot q^i) \in TQ$, and the Lagrangian $L : TQ \to \mathbb{R}$, but we are not given the constraints in Definition \ref{def:nhc-constraint}. The goal is to find the constraint distribution, or the collection of all horizontal space $\mathcal{D} = \sqcup_q \, H_q$ (explained in Section \ref{sec:connections}) of the system, and the full dynamics of the system.

\section{Background}
We will discuss some mathematical preliminaries for nonholonomic systems, symmetry, Lagrangian reduction, and the Hamel's formalism for mechanical systems in a moving basis.

\subsection{Connections and the Horizontal Lift}
\label{sec:connections}

A connection is a mathematical object that helps one understand a nonholonomic constraint as a map between spaces, which leads to the following discussion of the distribution under constraints, and the horizontal lift as an action of mapping any velocity to the constraint distribution. We introduce the concept of fiber bundle first to give a general geometric structure for nonholonomic systems.

\begin{definition}[\cite{Lee2000} Page 268]
\textit{
    Let $Q_s$ and $Q_r$ be topological spaces. A \textbf{fiber bundle} over $Q_s$ with model fiber $Q_r$ is a topological space $Q$ with a surjective continuous map $\pi : Q \to Q_s$ with the property that for each $s \in Q_s$, there exist a neighborhood $U$ of $s$ in $Q_s$ and a homeomorphism $\Phi : \pi^{-1}(U) \to U \times Q_r$, called a local trivialization of $Q$ over $U$.
}
\end{definition}

Considering the tangent map $T_q\pi$, if we split $T_qQ$ into separate parts $T_qQ = H_q \times V_q$, then with $T_q\pi : T_qQ \to H_q$, $TQ$ is a fiber bundle, where the fiber $V_q$ lies in the kernel of the map $T_q\pi$. Usually we denote $H_q$ the \textbf{horizontal space}, $V_q$ the \textbf{vertical space}, and $T_q\pi$ the \textbf{horizontal lift} $\mathrm{hor}$.

\begin{definition}[\cite{Bloch2015}]
    An \textbf{Ehresmann connection} A is a vector-valued one-form on $Q$ that satisfies:
    (1) $A$ is vertical valued: $A_q : T_qQ \to V_q$ is a linear map for each point $q \in Q$.
    (2) $A$ is a projection: $A(v_q) = v_q$ for all $v_q \in V_q$.
\end{definition}

One can see that $H_q$ is the kernel of $A_q$. Suppose the configuration space $Q$, or $Q_r \times Q_s$ has coordinates $(r^\alpha, s^a)$ \footnote{In this paper we assume Einstein summation convention for simplicity.}. Then we represent the connection in the system as
\begin{equation}
    A = \omega^{a}\frac{\partial}{\partial s^{a}}, \, \mathrm{where} \,\, \omega^{a}(q) = ds^{a} + \mathcal{A}^{a}_{\alpha}(r,s)dr^{\alpha},
\end{equation}
Then for a given $v_q = \dot{r}^{\alpha}\frac{\partial}{\partial r^{\alpha}} + \dot{s}^{a}\frac{\partial}{\partial s^{a}}$ in $T_qQ$, 
\begin{equation}
\label{eq:A(vq)}
    A_q(v_q) = (\dot{s}^{a} + \mathcal{A}^{a}_{\alpha} \dot{r}^{\alpha})\frac{\partial}{\partial s^{a}}.
\end{equation}
We can see that $A$ is a tangent projection map and gives a fiber bundle $T_qQ$. Correspondingly, the \textbf{horizontal lift} can be defined by subtracting the vector $v_q$ by the vertical component,
\begin{equation}
\label{eq:hor}
    \mathrm{hor}\,v_q = v_q - A_q(v_q) = \dot{r}^{\alpha} \frac{\partial}{\partial r^{\alpha}} - \mathcal{A}_{\alpha}^{a}\dot{r}^{\alpha} \frac{\partial}{\partial s^{a}}.
\end{equation}

\subsection{Matrix Representation of Connections and the Horizontal Lift}
\label{sec:matrix}
One way to represent the connection $A$ is through matrix operators\footnote{This is also what we write in computer programs.}, by considering each differential form $\omega^a$ as an operator in the row space, mapping a vector in the column space to a real value, and then the value is assigned to $\frac{\partial}{\partial s^a}$. More specifically, if we let $s$ and $r$ also represent the dimensionality of $Q_s$ and $Q_r$,
\begin{equation}
    A = \begin{bmatrix}
        I^{s\times s} & \mathcal{A}^{s\times r} \\
        \mathbf{0}^{r\times s} & \mathbf{0}^{r\times r}
    \end{bmatrix},
\end{equation}
where the previous notation $\mathcal{A}^i_j$ in Section \ref{sec:connections} simply represents the element of $\mathcal{A}$ at the $i^\mathrm{th}$ row and $j^\mathrm{th}$ column, and $A$ is the \textbf{reduced row echelon form} of matrix $a$ in Definition \ref{def:nhc-constraint}. For a random velocity $v_q = \begin{bmatrix}\dot{s} & \dot{r}\end{bmatrix}^T$ in $T_qQ$,
\begin{equation}
    A_q(v_q) = \begin{bmatrix}
        \dot{s} + \mathcal{A}(r,s)\dot{r} \\
        \mathbf{0}^{r\times 1}
    \end{bmatrix} \in V_q \,.
\end{equation}
Similarly, we can represent the horizontal map by subtracting $A$ from $I$,
\begin{equation}
    \mathrm{hor} = I - A = \begin{bmatrix}
        \mathbf{0}^{s\times s} & -\mathcal{A}^{s\times r} \\
        \mathbf{0}^{r\times s} & I^{r\times r}
    \end{bmatrix},
\end{equation}
which can act on $v_q$ and result in
\begin{equation}
    \mathrm{hor}_q\, v_q = \begin{bmatrix}
        -\mathcal{A}(r,s)\dot{r} \\
        \dot{r}
    \end{bmatrix} = \begin{bmatrix}
        -\mathcal{A}(r,s) \\ 
        I^{r\times r}
    \end{bmatrix} \dot{r} \in H_q \,.
\end{equation}
Since $\dot{r}$ represents any vector of $r$ dimension, we have
\begin{equation}
\label{eq:Hq}
    H_q = \mathrm{Col}\Big(\begin{bmatrix}
    -\mathcal{A}(r,s) \\ 
    I^{r\times r}
    \end{bmatrix}\Big).
\end{equation}

\begin{proposition}
\label{prop:Ahor}
    (1) Both $A$ and $hor$ are projections. (2) They lie in the kernel of each other. (3) Both $V_q$ and $H_q$ are unique given a nonholonomic constraint \eqref{eq:def-constraint}.
    \begin{proof}
        It is easy to check that (1) $A^2 = A$ and $\mathrm{hor}\circ\mathrm{hor}=\mathrm{hor}$, and (2) $A\,\mathrm{hor}=\mathrm{hor}\,A=0$. For (3), consider that for any given constraint $a(q)$, its reduced row echelon form $A_q$ is unique, so $\mathrm{hor}$ is unique. Therefore, $V_q$ and $H_q$ are unique.
\end{proof}
\end{proposition}

\subsection{Reduction to Lie Algebra}
\label{sec:reduction}
Now suppose the system is partially evolving on a Lie group (refer to Appendix \ref{app:vector} for a more rigorous definition of vector fields, Lie groups, Lie algebras), and we can split the configuration space to $Q/G \times G$ with points $(r, g)$ in contrast to $Q_r \times Q_s$ with points $(r, s)$. Suppose that the basis $e_a \in \mathfrak{g}, \, a=1,\dots,k$ spans the entire Lie algebra, so all the vectors on $G$ generated by left translation $L_g^*e_a, \, g\in G$ form a set of left-invariant Lie algebra vector fields. Rewriting the horizontal space \eqref{eq:Hq} in the coordinates corresponding to  $\frac{\partial}{\partial r^\alpha}$ and $L_g^*e_a$, we get
\begin{equation}
\label{eq:u_alpha}
    u_{\alpha} = \frac{\partial}{\partial r^\alpha} - \mathcal{A}(r, g)L_g^*e_a, \quad \alpha = 1, \dots, \sigma.
\end{equation}
which can be seen as a new set of coordinates that span either partially or the entire constraint distribution $\mathcal{D}$, and $\mathcal{A}$ becomes a $\mathfrak{g}$-valued function. For the part of vector field that is $G$-invariant, we can still represent it by  the vector fields
\begin{equation}
\label{eq:u_sigma}
    u_{\sigma+a} = L_g^*e_a, \quad a = 1, \dots, k.
\end{equation}
With $u_i$ for $i=1,\dots,\sigma+k$, we can represent any vector on $TQ/G \times \mathfrak{g}$. Therefore, these vector fields are sufficient to represent the entire vector field of a nonholonomic system with connection $A$ and with a Lagrangian that is invariant under the group action $L_g$. 

The components of the coordinate $(u_\alpha, u_{\sigma+a})$ can be represented as $(\dot{r}^\alpha, \Omega^a)$, where $\Omega^a = \xi^a + \mathcal{A}^a_\alpha \dot{r^\alpha}$ and $\xi = L_{g^{-1}}^*\dot{g}$ are body velocities with and without connections. Here we can understand $(\dot{r}^\alpha, \Omega^a)$ as the horizontal and vertical velocities (recall \eqref{eq:A(vq)}, \eqref{eq:hor}). The $G$-invariant Lagrangian $L(r, \dot{r}, g, \dot{g})$ can be reduced to $l(r, \dot{r}, \Omega)$ or $l(r, \dot{r}, \xi)$ by pulling back the global velocities in $TG$ to the body velocities in $\mathfrak{g}$. Furthermore, the system must satisfy the Hamel equations:

\subsection{The Hamel Equations}
The original Hamel equations were derived in a simpler setting, in particular to deal with the Euler-Lagrange equations in a moving basis. While the above description of the reduced nonholonomic systems can be interpreted in the moving basis \eqref{eq:u_alpha} and \eqref{eq:u_sigma}, there are some differences in coordinate selections compared to the original Hamel equations. For readability of this paper, we introduce the original Hamel equations and let the reader refer to \cite{Bloch2009Quasi} for a more detailed derivation for the nonholonomic version.

Let $v = (v^1,\dots, v^n) \in \mathbb{R}$ be components of the velocity vector $\dot{q} \in TQ$ relative to the basis $u_1,\dots,u_n$, i.e., $\dot{q} = v^i u_i(q)$, so that the Lagrangian $L(q,\dot{q}) = L(q, v^iu_i(q))$ can be reduced to $l(q,v)$. 
\begin{theorem}[Hamel Equations]
\label{thm:hamel}
    The evolution of the variables $(q, v)$ satisfying Hamilton's principle is governed by the Hamel equations 
    \begin{equation}
    \label{eq:hamel}
        \frac{\mathrm{d}}{\mathrm{d}t}\frac{\partial l}{\partial v} = \Big[v, \frac{\partial l}{\partial v}\Big]_q^* + u[l]
    \end{equation}
    In the coordinate form, 
    \begin{equation}
        \frac{\mathrm{d}}{\mathrm{d}t}\frac{\partial l}{\partial v^i} = c_{ji}^m\frac{\partial l}{\partial v^m}v^j + u_i[l],
    \end{equation}
    where the Lie algebra structure constants $c_{ji}^m$ are defined by
    \begin{equation}
        [u_j, u_i] = c_{ji}^m u_m.
    \end{equation}
\end{theorem}

The velocities $v^i$ in our case are $(\dot{r}, \xi^a)$ or $(\dot{r}, \Omega^a)$. While it is non-trivial to obtain the same formalism for nonholonomic systems with symmetry, one can in general apply the coordinates defined in \eqref{eq:u_alpha} and \eqref{eq:u_sigma} to Theorem \ref{thm:hamel} to generate the full dynamics of momentum $(\frac{\partial l}{\partial \dot{r}}, \frac{\partial l}{\partial \xi})$ or $(\frac{\partial l}{\partial \dot{r}}, \frac{\partial l}{\partial \Omega})$. In the rest of the paper, we also denote them $p_{\dot{r}}$, $p_\xi$ and $p_\Omega$. Notice that the difference between the dynamics of $(p_{\dot{r}}, p_\xi)$ and $(p_{\dot{r}}, p_\Omega)$ is that the former  gives a reduced system with symmetry while  the latter  gives a reduced system with symmetry and also with the connection.

\section{Learning Nonholonomic Dynamics}
\label{sec:learning-nhc}
Here we state our main contribution of this paper: a general procedure to learn nonholonomic dynamics with known Lagrangian $L(q, \dot{q})$ and unknown constraint $a$. From now on, we represent any element parameterized by a neural network with a tilde symbol $\tilde{\cdot}$.

\subsection{Step One: Data Reduction to Lie Algebra}
For a general nonholonomic systems on $TQ$ with data set $(q, \dot{q})^i$, we first split the system into the product of a known Lie group $G$ and the quotient of $TQ$ by this group $TQ/G$, so the data is represented as $(r, \dot{r}, g, \dot{g})^i$. For each data point, we can pull back the group elements by $\xi^i = L_{(g^i)^{-1}}^*g^i$, and the data set becomes $(r, \dot{r}, \xi)^i$.

\subsection{Step Two: Generate Moving Basis and Quasi-velocities}
Since we already chose the Lie group $G$, we can generate the Lie algebra vector field by push forward of  a frame $e_a$ defined in Section \ref{sec:reduction}, $L_g^*e_a$. This also corresponds to the vector fields $u_{\sigma+a}$ defined in \eqref{eq:u_sigma}. It can be seen that $u_{\sigma+a}$ is known once a group $G$ is chosen.

Now to have $u_\alpha$ defined in \eqref{eq:u_alpha}, we parameterize $\mathcal{A}(r,g)$ with a neural network, so that
\begin{equation}
    \tilde{u}_\alpha = \frac{\partial}{\partial r^\alpha} - \tilde{\mathcal{A}}(r,g)L_g^*e_a.
\end{equation}
Once we have the parameterized $\tilde{u}_\alpha$ and $u_{\sigma+a}$, we can compute the parameterized quasi-velocities $\tilde{\Omega}^a = \xi^a + \tilde{\mathcal{A}}^a_\alpha\dot{r}^\alpha$.

\subsection{Step Three: Compute Momentum and Generate Dynamics}
Knowing the Lagrangian $L(q, \dot{q})$, we can first reduce it to $l(r,\dot{r}, \xi) = L_{g^{-1}}^*L(q,\dot{q})$. We check by the  chain rule that $\frac{\partial l}{\partial \xi} = \frac{\partial l}{\partial (\xi + \tilde{\mathcal{A}}\dot{r})} \frac{\partial (\xi + \tilde{\mathcal{A}}\dot{r})}{\partial \xi} = \frac{\partial l}{\partial \tilde\Omega}$. 
So the momentum $\{p_\Omega^i\}$ can be directly calculated from the known $l(r, \dot{r}, \xi)$ and $\xi^i$. The momentum $\{ p_{\dot{r}}^i \}$ can be calculated by the chain rule
\begin{equation}
\label{eq:p_rdot_tilda}
    \tilde p_{\dot{r}} = \frac{\partial l}{\partial \dot{r}} = \frac{\partial l_{\dot{r}}}{\partial \dot{r}} + \big<\frac{\partial l}{\partial \xi}, - \tilde{\mathcal{A}} \big>,
\end{equation}
where $l_{\dot{r}}$ represents the part of function in $l$ that does not contain $\xi$. An interesting observation here is that if one takes the time derivative of the second part of the equation $\big<\frac{\partial l}{\partial \xi}, - \tilde{\mathcal{A}} \big>$, it is equivalent to the nonholonomic momentum equation. However, this will be beyond the scope of this paper's discussion.

\subsection{Step Four: Training with Neural ODE}
With the data points $\{p^i\} = \{ (p_{\dot{r}}, p_\Omega)^i \}$, and the dynamics calculated by substituting $u_\alpha$, $u_{\sigma+a}$ into Theorem \ref{thm:hamel}, we can train the dynamics as a neural ODE, with loss
\begin{equation}
\label{eq:loss}
    \mathcal{L} = \sum_i \Big\| \int_{t_i}^{t_{i+1}} \frac{\mathrm{d}}{\mathrm{d}t} \tilde p^i \mathrm{d}t - p^{i+1} \Big\|_2^2.
\end{equation}
In this formulation, the dynamics is generated through parameterizing only the nonholonomic constraint, and the constraint is fully represented in an unconstrained ODE. After training, the constraint is therefore implicitly discovered as the output of the network $\tilde A(r, g)$. To guarantee local convergence of this training process, we prove the following theorem.

\begin{theorem}[Unique Dynamics Generated by $\mathcal{A}$]
\label{thm:uniqueA}
    Consider a nonholonomic system with left invariant Lagrangian $L(q, \dot{q})$, and nonholonomic constraint $a(q)$, evolving on  $TQ/G \times G$. Define its coordinates as in \eqref{eq:u_alpha} and \eqref{eq:u_sigma}, then its dynamics generated by \eqref{eq:hamel} is defined locally uniquely by the connection.
\end{theorem}
We leave the proof of the theorem in Appendix \ref{app:proof} for readers who are interested. This theorem guarantees that there is a minimum for the loss function \eqref{eq:loss}.

\section{Learning the Rolling Disk Dynamics}
\label{sec:rolling-penny}

In this section, we discuss our approach applied to the rolling disk dynamics \cite{Bloch2015}. The rolling disk dynamics, by definition, is the dynamics of a penny or a disk-like object rolling on a flat surface without slipping. In our paper, we also assume the penny does not fall. The configuration space is denoted by $(\theta, \varphi, x, y, \dot{\theta}, \dot{\varphi}, \dot{x}, \dot{y})$, where $\theta$ is the angle of one fixed point on the penny relative to the horizontal line, $\varphi$ is the rolling direction of the penny relative to the horizontal coordinate of the plane, and $x$, $y$ are the location of the penny on the plane. The Lagrangian is given by the kinetic energy
\begin{equation}
    L(q, \dot{q}) = \frac{1}{2} I \dot{\theta}^2 + \frac{1}{2} J \dot{\varphi}^2 + \frac{1}{2} m (\dot{x}^2 + \dot{y}^2),
\end{equation}
where $I$ is the inertia of the penny about is center of rotation, $J$ is the inertia of the penny about $\varphi$, and $m$ is the mass of the penny. The non-slip condition gives the constraint \eqref{eq:constraint}, which is non-integrable.

\subsection{The $\mathbb{R}\times SE(2)$ version}
\label{sec:RxSE2}

For demonstration purposes, we derive the learning process of the system on the configuration space $\mathbb{R} \times SE(2)$, namely $r = \theta$ and $g = (\varphi, x, y)$. The group action of $SE(2)$ gives the left-invariant Lie algebra vector field
\begin{align}
    u_2 &= \frac{\partial}{\partial \varphi} - y \frac{\partial}{\partial x} + x \frac{\partial}{\partial y}, \\
    u_3 &= \frac{\partial}{\partial x}, \\
    u_4 &= \frac{\partial}{\partial y},
\end{align}
corresponding to \eqref{eq:u_sigma}. So for a given velocity at time $t_i$, $(g, \dot{g})^i$, we compute the Lie algebra element as $(\xi^1, \xi^2, \xi^3)^i = (\dot{\varphi}, \dot{x} + y\dot{\varphi}, \dot{y} - x\dot{\varphi})^i$.

Now we come to Step Two, where we parameterize $\mathcal{A}(r, g)$ and define $u_1$.
\begin{equation}
    \tilde u_1 = \frac{\partial}{\partial r} - \tilde{\mathcal{A}}^a(r, g) u_{1+a}, \quad a = 1, 2, 3.
\end{equation}
Together with $(\xi_1, \xi_2, \xi_3)$, we can compute the parameterized quasi velocities $\Omega^a = \xi^a + \tilde{\mathcal{A}}^a \dot{\theta}$.

Then we compute the dynamics as in Step Three. Substituting $\xi^a$ to the Lagrangian, we get the reduced Lagrangian,
\begin{equation}
    l(\theta, \dot{\theta}, \xi^1, \xi^2, \xi^3) = \frac{1}{2}I \dot{\theta}^2 + \frac{1}{2} J (\xi^1)^2 + \frac{1}{2} m \big((-y\xi^1 + \xi^2)^2 + (x\xi^1 + \xi^3)^2\big),
\end{equation}
so we can calculate the momentum by
\begin{align}
\label{eq:p}
    p_2 &= \frac{\partial l}{\partial \Omega^1} = \frac{\partial l}{\partial \xi^1} = J\xi^1 - my(-y\xi^1 + \xi^2) + mx(x\xi^2 + \xi^3), \\
    p_3 &= \frac{\partial l}{\partial \Omega^2} = \frac{\partial l}{\partial \xi^2} = m(-y\xi^1 + \xi^2), \\
    p_4 &= \frac{\partial l}{\partial \Omega^3} = \frac{\partial l}{\partial \xi^3} = m(x\xi^1 + \xi^3), \\
    p_1 &= \frac{\partial l}{\partial \dot{\theta}} = I \dot{\theta}.
\end{align}
Note that these are the ground truth momenta. In the parameterized momenta, $p_2, p_3, p_4$ stay the same, while we use \eqref{eq:p_rdot_tilda} to compute $p_1$. Then we compute the Lie bracket coefficients with the coordinates $u_i$.

\begin{align}
    [u_2, u_1] &= - (\frac{\partial \tilde{\mathcal{A}^a}}{\partial \varphi} - y \frac{\partial \tilde{\mathcal{A}^a}}{\partial x} + x \frac{\partial \tilde{\mathcal{A}^a}}{\partial y}) u_{1+a} + \mathcal{A}^a u_{1+a} [u_2] = \tilde{c}_{21}^m u_m, \\
    [u_2, u_3] &= - \frac{\partial}{\partial y} = -1 \cdot u_4 = c_{23}^4 u_4, \\
    [u_2, u_4] &= \frac{\partial}{\partial x} = 1 \cdot u_3 = c_{24}^3 u_3, \\
    &\dots \nonumber
\end{align}
We omit the entire calculation for simplicity. The general idea is to calculate $c_{ji}^m$, where some of them are parameterized and some of them are not. Then we generate the dynamics based on \eqref{eq:hamel}, where we see $v^1$ as $\dot{r}$ and $v^2, v^3, v^4$ as $\Omega^1, \Omega^2, \Omega^3$. In Appendix \ref{app:rolling-penny}, we derive the dynamics with known constraints to familiarize the readers with the Hamel equations.

It is worth mentioning that this approach is agnostic to your choice of group $G$. That is, one will still recover the momentum dynamics by choosing the configuration $\mathbb{R} \times S^1 \times \mathbb{R}^2$ with coordinate $(\varphi, \dot{\varphi}, g, \dot{g})$ and $g = (\theta, x, y)$. However, you are evaluating a different momentum about a different moving basis, and the map $\mathcal{A}$ will appear to be different. By changing coordinates, the system will still follow the same trajectory.

\subsection{Results}

The neural network chosen for this task is a Feed-forward Neural Network\footnote{Code can be found on the \href{https://github.com/Wangbaiyue007/learning-rolling-penny/tree/hamel}{GitHub} page.} with 3 hidden layers and a hidden dimension of 20. Each hidden layer has a linear layer and a nonlinear layer of the \textit{Sine()} function, and the last nonlinear layer is chosen to be the \textit{ELU()} function that allows near zero output. The training set is calculated as described in Section \ref{sec:rolling-penny} with random system parameters $I, J, m, R$ and initial conditions, in the time ranging from $0$ to $20 \, sec$ with a step size of $0.01 \, sec$. The entire training process involves $2000$ epochs. The program was carried out with \textit{PyTorch} and \textit{torchdiffeq} \cite{chen2019neuralode}.

\subsubsection{Learning on $\mathbb{R} \times SE(2)$}

\begin{figure}[t]
    \scalebox{.5}{\input{figs_2026/2000.pgf}}
    \caption{The above figure shows the training result after 2000 epochs. The sub plot on the left shows the ground truth momenta (solid lines) and the learned dynamics (dashed lines) generated by a numerical integrator. The sub plot on the right shows the ground truth $\mathcal{A}(r,g)$ (solid lines) and the parameterized $\tilde{\mathcal{A}}(r,g)$ (dashed lines) evaluated after training.}
    \label{fig:2000}
\end{figure}

The results of our training process carried out on $\mathbb{R} \times SE(2)$ are shown in Figure \ref{fig:2000}. The ground truth $\mathcal{A}^a(r, g)$ is defined by $\mathcal{A}^1(\theta, \varphi, x, y) = 0$, $\mathcal{A}^2(\theta, \varphi, x, y) = -R\cos{\varphi}$, and $\mathcal{A}^3(\theta, \varphi, x, y) = R\sin{\varphi}$. We can see that the learned momentum converges to the ground truth, and the learned constraint $\tilde{\mathcal{A}}^a(r,g)$ converges to $\mathcal{A}^a(r,g)$ by evaluation. Notice that these constraints were learned completely from the data set $\{(r, \dot r, g, \dot{g})^i\}$ without knowing the ground truth $\mathcal{A}^a(r,g)$. Therefore, we call this result "constraint discovery". Any base line methods that parameterize the dynamics directly do not imply any information about the constraint or the symmetry of the system.

Notice that the constraint does not  perfectly match the ground truth (this is further discussed in our conclusions).
We think that this is due to the following: the numerical integration method introduces some error; the training process only includes one trajectory with one initial condition, so in some cases the momentum trajectory does not indicate enough about the constraint. Actually, although the entire vector field is unique according to Theorem \ref{thm:uniqueA}, a single trajectory can be generated by different vector fields. However, based on the uniqueness of the vector field, the convergence can be improved to reach a sufficiently small error by training on different trajectories.

\section{Conclusions and Future Work}
\label{sec:conclusion}

In this paper we proposed a general learning approach for learning the dynamics of a nonholonomic system subject to unknown constraints, where the results imply that the parameterized constraint is implicitly learned by training the system in a neural ODE fashion represented by Hamel equations. This is one of the first approaches to learning a nonholonomic system with symmetry and there are some open problems that the authors seek to improve in future work.

\noindent\textbf{Impact of Training Set} The proposed method in this paper aims at learning the constraint based on a single trajectory of many points starting from the same initial condition. Therefore, a badly chosen initial condition may lead to results that are not perfectly converging. We need to find a way to learn from different initial conditions, without relying on a single trajectory, to address this problem. 
Furthermore, we find it interesting to study the global convergence of the network.

\noindent\textbf{Numerical Integration error} Since the entire dynamics of the system is calculated by numerical ODE integration, the error caused by the integrator cannot be avoided. We expect to improve future results by considering discrete nonholonomic dynamics \cite{cortes2001, colombo2015} with a different integration approach that preserves symmetry properties, or the momentum equation. There are some approaches for discrete Lagrangian systems \cite{hansen2025learning}.

\noindent\textbf{Application to Real-world Systems} nonholonomic systems exist broadly in the real world, including robotics, vehicle dynamics, control theory, aerospace, marine, and computer graphics. Having some observed data from those systems and discovering the constraint would be an interesting application of the proposed learning process.

\acks{Partially supported by NSF grant  DMS-2103026, and AFOSR grants FA
9550-22-1-0215 and FA 9550-23-1-0400.}

\bibliography{references}

\newpage
\appendix
\section{Vector Fields, Lie Groups, Lie Algebras}
\label{app:vector}
    Let $Q$ be a smooth n-manifold and let $q\in Q$. Then $T_qQ$ is an n-dimensional vector space. Let the coordinate vectors $\frac{\partial}{\partial q^1}, \dots, \frac{\partial}{\partial q^n}$ form a basis for $T_qQ$. A tangent vector $v \in T_qQ$ can be written uniquely as a linear combination (see \cite{Lee2000} page 61)
    \begin{equation}
    \label{eq:v}
        v = v^i \frac{\partial}{\partial q^i} \Bigg\vert_q,
    \end{equation}
    where the coeffcients  $(v^1, \dots , v^n)$ are called the components of $v$. A vector field on $Q$ is a continuous map $u: Q \to TQ$ with the property that $u(q) \in T_qQ$ for each $q\in Q$. Let $u_i$, $i=1,\dots, n$ be smooth independent local vector fields defined by (see \cite{Bloch2009Quasi})
    \begin{equation}
        u_i(q) = \psi^j_i (q) \frac{\partial}{\partial q^j}, \quad i,j = 1,\dots, n.
    \end{equation}
    Then the $u_i$ form a basis for $TQ$. With the $u_i$ as a basis, and again using $v^i$ as components as in \eqref{eq:v}, a velocity vector $\dot{q} \in T_qQ$ can be written as
    \begin{equation}        \dot{q} = v^i u_i(q),
    \end{equation}
    then for a function $L: TQ \to \mathbb{R}$, we can define the reduced Lagrangian $l: TQ \to \mathbb{R}$
    \begin{equation}
        l(q, v) := L(q, v^iu_i(q)).
    \end{equation}
    We can view the  vector fields $u_i$  as operators and define the directional derivatives $u_i[l]$ as
    \begin{equation}
        u_i[l] = \psi^j_i\frac{\partial l}{\partial q^j}.
    \end{equation}
    The evolution of the variables $(q, v)$ satisfying Hamilton's principle is governed by the Hamel equations \eqref{eq:hamel} (\cite{Bloch2009Quasi}).

    Now we consider systems on Lie groups. A Lie group is a smooth manifold $G$ that is also a group in the algebraic sense (\cite{Lee2000} Chapter 7). If $G$ is a Lie group, any element $g \in G$ defines maps $L_g: G \to G$ called left translation and $L_g^*: TG \to TG$ called pushforward,
    \begin{align}
        L_g(h) = gh, &\quad h \in G, \\
        L_g^*(v_h) = (dL_g)_h(v_h), &\quad v_h \in T_hG.
    \end{align}
    Let $X$, $Y$ be smooth vector fields on $Q$, and $X = X^i \frac{\partial}{\partial q^i}$, $Y = Y^i \frac{\partial}{\partial q^j}$. Then a Lie bracket is defined by (\cite{Lee2000} Proposition 8.26)
    \begin{equation}
        [X, Y] = \Big( X^i \frac{\partial Y^j}{\partial q^i} - Y^i \frac{\partial X^j}{\partial q^i} \Big) \frac{\partial}{\partial x^j}.
    \end{equation}
    A Lie algebra is a real vector space $\mathfrak{g}$ with a Lie bracket $\mathfrak{g} \times \mathfrak{g} \to \mathfrak{g}$ that satisfies bilinearity, antisymmetry, and the Jacobi identity (\cite{Lee2000} Page 190). The Lie algebra of all smooth left-invariant vector field is called the Lie algebra of $G$, which is also isomorphic to $T_eG$.

\section{Proof of Theorem \ref{thm:uniqueA}}
\label{app:proof}
    To prove this theorem, it suffices to say that under some small perturbation of $\mathcal{A}(r,g)$, the generated dynamics is different. To make the proof cleaner, we proceed in a more general setting on the configuration space $Q$. Define a basis for $Q$
    \begin{align}
        u_\alpha(q) &= \psi^j_\alpha(q)\frac{\partial}{\partial q^j}, \quad j = 1, \dots, n, \quad \alpha = 1,\dots,\sigma \\
        u_{\sigma+a}(q) &= \psi^j_{\sigma+a}(q) \frac{\partial}{\partial q^j},\quad j = m,\dots,n \quad a = 1, \dots,n-\sigma
    \end{align}
    and the velocity is
    \begin{equation}
        \dot{q} = v^\alpha u_\alpha + v^{\sigma + a}u_{\sigma + a}.
    \end{equation}
    Assume that the reduced Lagrangian is $l(q, v)$ and the nominal dynamics is generated by
    \begin{equation}
        \frac{\mathrm{d}}{\mathrm{d}t} \frac{\partial l}{\partial v^i} = c_{ji}^m \frac{\partial l}{\partial v^m} v^j + u_i[l].
    \end{equation}
    Then, we assume that if we perturb the coordinates by a function $\delta$ of $q$
    \begin{equation}
        u_\alpha'(q) = \psi^j_\alpha \frac{\partial}{\partial q^j} - \delta^{\sigma+a}_\alpha u_{\sigma+a},
    \end{equation}
    the generated dynamics stays the same. The velocity after this modification is
    \begin{equation}
        \dot{q} = v^\alpha u_\alpha' + w^{\sigma+a} u_{\sigma+a} := v^\alpha u_\alpha' + (v^{\sigma+a} + \delta^{\sigma+a}) u_{\sigma+a},
    \end{equation}
    where $\delta^{\sigma+a} = \sum_\alpha \delta^{\sigma+a}_\alpha$. By the  chain rule, $\frac{\partial l}{\partial v^{\sigma+a}} = \frac{\partial l}{\partial w^{\sigma+a}}$. Now considering the case $i = \alpha, j = \sigma+a$ in Hamel equations,
    \begin{equation}
        \frac{\mathrm{d}}{\mathrm{d}t} \frac{\partial l}{\partial v^\alpha} = c_{\sigma+a, \alpha}'^m \frac{\partial l}{\partial v^m} w^{\sigma+a} + c_{\beta, \alpha}'^m \frac{\partial l}{\partial v^m} v^{\beta} + u_\alpha'[l],
    \end{equation}
    one can see that if $u_\alpha [l] = u_\alpha'[l]$, then $\delta^{\sigma+a}_\alpha u_{\sigma+a}[l] = 0$. So if $l$ is a function of $q^i, i = m,\dots,n$, then $\delta^{\sigma+a}_\alpha = 0$. In the case where $l$ is not a function of $q^i, i = m,\dots,n$, one needs to check the rest of the equation. Now we check $c_{\sigma+a, \alpha}'^m$
    \begin{align}
        [u_{\sigma+a}, u_\alpha '] &= [u_{\sigma+a}, u_\alpha - \delta^{\sigma+a}_\alpha u_{\sigma+a}] \\
        &= [u_{\sigma+a}, u_\alpha] + [u_{\sigma+a}, -\delta^{\sigma+a}_\alpha u_{\sigma+a}] \\
        &:= c_{\sigma+a, \alpha}^m u_m + b_{\sigma+a,\alpha}u_{\sigma+a};
    \end{align}
    compared to the nominal case, there is an extra term $b_{\sigma+a,\alpha}u_{\sigma+a}$. We examine this term
    \begin{align}
        [u_{\sigma+a}, -\delta^{\sigma+a}_\alpha u_{\sigma+a}] &= -\psi^j_{\sigma+a}\frac{\partial \delta^{\sigma+a}\psi^l_{\sigma+a}}{\partial q^j} \frac{\partial}{\partial q^l} + \delta^{\sigma+a}_\alpha\psi^j_{\sigma+a} \frac{\partial \psi^l_{\sigma+a}}{\partial q^j} \frac{\partial}{\partial q^l} \\
        &= -\psi^j_{\sigma+a} \big(\frac{\partial\delta^{\sigma+a}_\alpha}{\partial q^j} \psi^l_{\sigma+a} + \delta^{\sigma+a}_\alpha\frac{\partial \psi^l_{\sigma+a}}{\partial q^j} \big)\frac{\partial}{\partial q^l} + \delta^{\sigma+a}_\alpha\psi^j_{\sigma+a} \frac{\partial \psi^l_{\sigma+a}}{\partial q^j} \frac{\partial}{\partial q^l} \\
        &= - \psi^j_{\sigma+a} \frac{\partial \delta^{\sigma+a}_\alpha}{\partial q^j} \psi^l_{\sigma+a} \frac{\partial}{\partial q^l} \\
        &= - \psi^j_{\sigma+a} \frac{\partial \delta^{\sigma+a}_\alpha}{\partial q^j} u_{\sigma+a} \\
        &= b_{\sigma+a,\alpha} u_{\sigma+a}.
    \end{align}
    For the term $c_{\sigma+a,\alpha}'^m p_m w^{\sigma+a}$ to equal the nominal $c_{\sigma+a,\alpha}^m p_m v^{\sigma+a}$, we subtract them and get
    \begin{equation}
        c_{\sigma+a,\alpha}'^m \frac{\partial l}{\partial v^m}\delta^{\sigma+a} + b_{\sigma+a,\alpha}\frac{\partial l}{\partial v^{\sigma+a}}(v^{\sigma+a} + \delta^{\sigma+a}) = 0,
    \end{equation}
    so for a small $\delta^{\sigma+a}$, 
    \begin{equation}
        b_{\sigma+a,\alpha}\frac{\partial l}{\partial v^{\sigma+a}}v^{\sigma+a} = 0,
    \end{equation}
    meaning
    \begin{align}
        b_{\sigma+a,\alpha} &= 0 \\
        \frac{\partial \delta^{\sigma+a}_\alpha}{\partial q^j} &= 0,
    \end{align}
    so $\delta^{\sigma+a}_\alpha$ is constant relative to $q$, which contradicts our assumption. Therefore, the perturbed coordinate must generate a different dynamics unless $\delta = 0$.

    As a general proof, this applies to our case in \eqref{eq:u_alpha}, \eqref{eq:u_sigma}. Together with Proposition \ref{prop:Ahor} (3), a nonholonomic system with constraint \eqref{eq:def-constraint} generates a locally unique dynamics.

\section{The Rolling Penny Dynamics with Known Constraint}
\label{app:rolling-penny}
    In Section \ref{sec:RxSE2}, we already derived the Lie algebra vector field $u_2, u_3, u_4$ and the parameterized horizontal vector field $\tilde{u}_1$. 
    
    With known constraint \eqref{eq:constraint}, we have
    \begin{equation}
        u_1 = \frac{\partial}{\partial \theta} + R\cos{\varphi}\frac{\partial}{\partial x} + R\sin{\varphi}\frac{\partial}{\partial y},
    \end{equation}
    so a vector that lies in the constraint distribution can be written
    \begin{equation}
        \dot{q} = \dot{\theta}u_1 + \Omega^1u_2 + \Omega^2u_3 + \Omega^3u_4,
    \end{equation}
    or
    \begin{equation}
        \dot{q} = \dot{\theta}\frac{\partial}{\partial \theta} + \xi^1u_2 + \xi^2u_3 + \xi^3u_4.
    \end{equation}
    Note that $\xi^1$ to $\xi^3$ can be calculated directly from data and from $u_2$ to $u_4$. Choose
    \begin{align}
        \Omega^1 &= \dot{\varphi} = \xi^1 + 0, \\
        \Omega^2 &= y\dot{\varphi} = \xi^2 + (-R\cos{\varphi}\dot{\theta}), \\
        \Omega^3 &= -x\dot{\varphi} = \xi^3 + (-R\sin{\varphi}\dot{\theta})
    \end{align}
    so that $\dot{\theta}u_1$ and $\Omega^1u_2+\Omega^2u_3+\Omega^3u_4$ span the constraint distribution 
    \begin{equation}
        \mathcal{D} = \big\{\frac{\partial}{\partial \theta} + R\cos{\varphi}\frac{\partial}{\partial x} + R\sin{\varphi}\frac{\partial}{\partial y}, \frac{\partial}{\partial \varphi} \big\}.
    \end{equation}

    Now we compute Lie brackets.
    \begin{align}
        [u_2, u_1] &= -R\sin{\varphi}\frac{\partial}{\partial x} + R\cos{\varphi}\frac{\partial}{\partial y} - R\cos{\varphi}\frac{\partial}{\partial y} + R\sin{\varphi}\frac{\partial}{\partial x} = 0, \\
        [u_3, u_1] &= 0, \\
        [u_4, u_1] &= 0, \\
        [u_3, u_2] &= \frac{\partial}{\partial y} = 1 \cdot u_4 = c_{32}^4 u_4, \\
        [u_4, u_2] &= - \frac{\partial}{\partial x} = -1 \cdot u_3 = c_{42}^3 u_3, \\
        [u_4, u_3] &= 0.
    \end{align}
    The anti-symmetric property of Lie brackets also gives $c_{ij}^m = - c_{ji}^m$. The values of $p_1$ to $p_4$ was given in \eqref{eq:p}. Now using the Hamel equations,
    \begin{align}
        \frac{\mathrm{d}}{\mathrm{d}t} p_1 &= 0, \\
        \frac{\mathrm{d}}{\mathrm{d}t} p_2 &= c_{32}^4p_4\Omega^2 + c_{42}^3p_3\Omega^3 = y\dot{\varphi}p_4 +x\dot{\varphi}p_3, \\
        \frac{\mathrm{d}}{\mathrm{d}t} p_3 &= c_{23}^4p_4\Omega^1 = -\dot{\varphi}p_4, \\
        \frac{\mathrm{d}}{\mathrm{d}t} p_4 &= c_{24}^3p_3\Omega^1 = \dot{\varphi}p_3,
    \end{align}
    and we obtain the dynamics of $p_1$ to $p_4$.

\end{document}